\documentclass[a4paper,11pt]{article}

\RequirePackage{amsmath,amsthm,amssymb}
\usepackage{hyperref}
\usepackage{setspace}

\newtheorem{thm}{Theorem}
\newtheorem{lem}[thm]{Lemma}
\newtheorem{rem}[thm]{Remark}

\newcommand{\I}{\textbf{I}}
\newcommand{\Ustirling}[2]{\genfrac{[}{]}{0pt}{}{#1}{#2}}
\newcommand{\forb}{\text{forb}}

\title{Inversion sequences avoiding the pattern 010}

\author{Benjamin Testart}

\usepackage[backend=bibtex]{biblatex}
\usepackage{a4wide}
\begin{document}
\renewcommand{\thefootnote}{\fnsymbol{footnote}}
\begin{center}
{\LARGE Inversion sequences avoiding the pattern 010} \\[15pt]
{\large Benjamin Testart\footnote[1]{Université de Lorraine, CNRS, Inria, LORIA, F 54000 Nancy, France}} \\
\end{center}
\renewcommand{\thefootnote}{\arabic{footnote}}
\abstract{Inversion sequences are integer sequences $(\sigma_1, \dots, \sigma_n)$ such that $0 \leqslant \sigma_i < i$ for all $1 \leqslant i \leqslant n$. The study of pattern-avoiding inversion sequences began in two independent articles by Mansour--Shattuck and Corteel--Martinez--Savage--Weselcouch in 2015 and 2016. These two initial articles solved the enumeration of inversion sequences avoiding a single pattern for every pattern of length 3 except the patterns 010 and 100. The case 100 was recently solved by Mansour and Yildirim. We solve the final case by making use of a decomposition of inversion sequences avoiding the pattern 010. Our decomposition needs to take into account the maximal value, and the number of distinct values occurring in the inversion sequence. We then expand our method to solve the enumeration of inversion sequences avoiding the pairs of patterns $\{010, 000\}, \{010, 110\}, \{010, 120\}$, and the Wilf-equivalent pairs $\{010, 201\} \sim \{010, 210\}$. For each family of pattern-avoiding inversion sequences considered, its enumeration requires the enumeration of some family of constrained words avoiding the same patterns, a question which we also solve.
} \medbreak

{\noindent \textbf{Keywords: }pattern avoidance, inversion sequences, enumeration, words, Stirling \break numbers}

\section{Introduction}

Let $\mathbb N$ be the set of natural numbers, including 0. Given a natural number $n \in \mathbb N$, we call \emph{integer sequences} of length $n$ the elements of $\mathbb N^n$. If $\sigma \in \mathbb N^n$ is an integer sequence, we write $\sigma = (\sigma_1, \dots, \sigma_n)$. We denote by $\I_n$ the set of \emph{inversion sequences} of length $n$, that is sequences $\sigma = (\sigma_1, \dots, \sigma_n)$ with elements in $\mathbb N$, such that $\sigma_i < i$ for all $i \in \{1, \dots, n\}$. There is a natural bijection between $\I_n$ and the set of permutations of $n$ elements, called the \emph{Lehmer code}, which explains the name "inversion sequence". If $\pi = (\pi_1, \dots, \pi_n)$ is a permutation, the inversion sequence $\sigma = (\sigma_1, \dots, \sigma_n)$ associated by the Lehmer code is defined by $\sigma_i = \#\{j \; | \; \pi_j > \pi_i \text{ and } j < i \}$ for all $i \in \{1, \dots, n\}$, i.e. $\sigma_i$ counts the number of inversions of $\pi$ whose second element is at position $i$.

The study of pattern-avoiding inversion sequences (and many more types of sequences avoiding patterns) branched from pattern-avoiding permutations, a well-established field of research in enumerative combinatorics, see \cite{Kitaev_2011}, \cite{Vatter_2015}, or \cite{wikipedia_permutations}.
Pattern-avoiding inversion sequences were first introduced in \cite{Mansour_Shattuck_2015} and \cite{Corteel_Martinez_Savage_Weselcouch_2016}. 
This study was extended using various definitions for 'patterns', see \cite{Martinez_Savage_2018} for patterns of relations, \cite{Yan_Lin_2020} for pairs of patterns, \cite{Auli_Elizalde_2019_1} for consecutive patterns, \cite{Auli_Elizalde_2019_2} for consecutive patterns of relations, \cite{Auli_Elizalde_2021} and \cite{Lin_Yan_2020} for vincular patterns, among other articles about the enumeration of pattern-avoiding inversion sequences.

In this paper, we will only consider 'classical' patterns: given two integer sequences $\sigma = (\sigma_1, \dots, \sigma_n) \in \mathbb N^n$ and $p = (p_1, \dots, p_k) \in \mathbb N^k$, we say that $\sigma$ \emph{contains} the pattern $p$ if and only if there exists a subsequence of $\sigma$ which is order-isomorphic to $p$. For example, the sequence $(4,3,2,5,4)$ contains the pattern 021 since both subsequences $(3,5,4)$ and $(2,5,4)$ are order-isomorphic to 021. A sequence \emph{avoids} a pattern $p$ if it does not contain $p$, e.g. the inversion sequence $(0,0,2,3,2,0,1,5)$ avoids the pattern 101.

\begin{table}[ht]
\begin{center}
    \begin{tabular}{|c|c|c|c|}
    \hline
    Pattern $p$ & $\#\I_n(p)$ for $n = 1, \dots, 7$ & Solved? & OEIS \cite{OEIS}\\
    \hline
    000 & 1, 2, 5, 16, 61, 272, 1385 & \cite{Corteel_Martinez_Savage_Weselcouch_2016} & A000111\\
    001 & 1, 2, 4, 8, 16, 32, 64 & \cite{Corteel_Martinez_Savage_Weselcouch_2016} & A000079\\
    011 & 1, 2, 5, 15, 52, 203, 877 & \cite{Corteel_Martinez_Savage_Weselcouch_2016} &  A000110\\
    012 & 1, 2, 5, 13, 34, 89, 233 & \cite{Corteel_Martinez_Savage_Weselcouch_2016} and \cite{Mansour_Shattuck_2015} & A001519\\
    021 & 1, 2, 6, 22, 90, 394, 1806 & \cite{Corteel_Martinez_Savage_Weselcouch_2016} and \cite{Mansour_Shattuck_2015} & A006318\\
    101 or 110 & 1, 2, 6, 23, 105, 549, 3207 & \cite{Corteel_Martinez_Savage_Weselcouch_2016} & A113227\\
    102 & 1, 2, 6, 22, 89, 381, 1694 & \cite{Mansour_Shattuck_2015} & A200753\\
    120 & 1, 2, 6, 23, 103, 515, 2803 & \cite{Mansour_Shattuck_2015} & A263778\\
    201 or 210 & 1, 2, 6, 24, 118, 674, 4306 & \cite{Mansour_Shattuck_2015} & A263777\\
    100 & 1, 2, 6, 23, 106, 565, 3399 & \cite{Mansour_Yildirim_2022} & A263780\\
    010 & 1, 2, 5, 15, 53, 215, 979 & Theorem \ref{thm010} & A263779\\
    \hline
    \end{tabular}
    \caption{Enumeration of inversion sequences avoiding a single pattern of length 3.}
    \label{table1}
\end{center}
\end{table}
If $P$ is a set of patterns, we denote by $\I_n(P)$ the set of inversion sequences of length $n$ avoiding all patterns in $P$. 
The enumeration of inversion sequences avoiding a single pattern of length 3 was already solved for all patterns except $010$ (see Table \ref{table1}). We solve the pattern 010 in Section \ref{010} by using a decomposition of inversion sequences avoiding 010, which involves constrained words avoiding 010. We first enumerate these pattern-avoiding words, exhibiting a new instance of the Stirling numbers. Then, our decomposition results in a recursive formula for inversion sequences avoiding 010.

\vspace{10pt}
\begin{table}[ht]
\begin{center}
    \begin{tabular}{|c|c|c|c|}
    \hline
    Pattern $p$ & $\#\I_n(010,p)$ for $n = 1, \dots, 7$ & Solved? & OEIS \cite{OEIS} \\
    \hline
    001 & 1, 2, 3, 4, 5, 6, 7 & \cite{Yan_Lin_2020} & A000027 \\
    011 & 1, 2, 4, 9, 23, 66, 210 & \cite{Yan_Lin_2020} & A026898 \\
    012 & 1, 2, 4, 8, 16, 32, 64 & \cite{Yan_Lin_2020} & A000079 \\
    021 & 1, 2, 5, 14, 42, 132, 429 & \cite{Yan_Lin_2020} & A000108 \\
    100 & 1, 2, 5, 15, 52, 203, 877 & \cite{Yan_Lin_2020} & A000110 \\
    101 & 1, 2, 5, 15, 52, 203, 877 & \cite{Martinez_Savage_2018} & A000110 \\
    000 & 1, 2, 4, 10, 29, 95, 345 & Theorem \ref{thm000} & A279552 \\
    120 & 1, 2, 5, 15, 52, 201, 845 & Theorem \ref{thm120} & A279559 \\
    201 or 210 & 1, 2, 5, 15, 53, 214, 958 & Theorem \ref{thm210} & A360052 \\
    110 & 1, 2, 5, 15, 52, 201, 847 & Theorem \ref{thm110} & A359191 \\
    102 & 1, 2, 5, 15, 51, 186, 707 & no & not yet \\
    \hline
    \end{tabular}
\caption{Enumeration of inversion sequences avoiding 010 and one other pattern of length 3.}
\label{table2}
\end{center}
\end{table}
The systematic study of inversion sequences avoiding pairs of patterns of length 3 was done by Lin and Yan in \cite{Yan_Lin_2020}. Among the pairs of patterns that contain 010, only five cases (up to Wilf-equivalence\footnote{We say that two patterns are \emph{Wilf-equivalent} if the corresponding avoidance classes of inversion sequences have the same enumeration sequence.}) are still unknown, see Table \ref{table2}. We solve four of those cases by adapting our approach for 010-avoiding inversion sequences, leaving only the pair $\{010, 102\}$ open.

\section{The pattern 010} \label{010}

The aim of this section is to enumerate inversion sequences avoiding the pattern 010. We solve this problem in Theorem \ref{thm010} by using a decomposition of 010-avoiding inversion sequences, in which a family of 010-avoiding words appears. We solve the enumeration of this family of words first, in Lemma \ref{lem010}.

Let $\Ustirling{n}{k}$ be the unsigned Stirling number of the first kind, counting permutations of length $n$ with $k$ cycles (among other combinatorial interpretations).
\begin{lem} \label{lem010}
Let $\mathfrak{A}_{n,k}$ be the set of 010-avoiding words $\omega = (\omega_1, \dots, \omega_n)$ of length $n$ on the alphabet $\{1, \dots, k\}$ such that $\omega$ contains all letters $\{1, \dots, k\}$ and $\omega_1 = k$. Let $\mathfrak{a}_{n,k} = \#\mathfrak{A}_{n,k}$. For all $n,k \geqslant 1$,
$$\mathfrak{a}_{n,k} = \Ustirling{n}{n+1-k}.$$
\end{lem}
\begin{proof}
    Let $n \geqslant k \geqslant 2$, and $\omega \in \mathfrak{A}_{n,k}$. Since $\omega$ avoids the pattern $010$, all letters 1 (which is the smallest letter) in $\omega$ are consecutive.
    \begin{itemize}
    \item If $\omega$ contains several letters 1, then removing one of them yields a word $\omega' \in \mathfrak{A}_{n-1,k}$, and this is clearly a bijection.
    \item If $\omega$ contains a single letter 1, then removing it and subtracting $1$ from all other letters yields a word $\omega' \in \mathfrak{A}_{n-1,k-1}$. Since $\omega_1 = k > 1$, there are $n-1$ possible positions for the letter 1 in $\omega$, so exactly $n-1$ words $\omega \in \mathfrak{A}_{n,k}$ yield the same $\omega'$.
    \end{itemize}
    Hence the following recurrence relation holds for all $n \geqslant k \geqslant 2$:
    $$\mathfrak{a}_{n,k} = \mathfrak{a}_{n-1,k} + (n-1) \cdot \mathfrak{a}_{n-1,k-1}.$$
    This recurrence relation is also satisfied by the Stirling numbers of the first kind $\Ustirling{n}{n+1-k}$, and we can easily verify that initial conditions also match.
\end{proof}
\begin{thm} \label{thm010}
Let $\mathfrak{B}_{n,m,d}$ be the set of 010-avoiding inversion sequences that have length $n$, maximum $m$, and contain exactly $d$ distinct values. Let $\mathfrak{b}_{n,m,d} = \#\mathfrak{B}_{n,m,d}$, so that we have
$$\# \emph{\I}_n(010) = \sum_{m = 0}^{n-1} \sum_{d=0}^n  \mathfrak{b}_{n,m,d}.$$
The numbers $\mathfrak{b}_{n,m,d}$ satisfy the following recurrence relation:
\begin{equation}
\mathfrak{b}_{n,m,d} =  \sum_{i = 0}^{d-1} \binom{m-i}{d-i-1} \sum_{p = m+1}^n \Ustirling{n-p+1}{n-p-d+i+2} \sum_{j = 0}^{m-1} \mathfrak{b}_{p-1,j,i}.
\label{eq010}
\end{equation}

\end{thm}
\begin{proof}
    Let $\sigma \in \mathfrak{B}_{n,m,d}$. Let $p$ be the left-most position of the value $m$ in $\sigma$. Let $\alpha = (\sigma_i)_{i \in \{1, \dots, p-1\}}$ be the subsequence of $\sigma$ to the left of the left-most $m$, and $\beta = (\sigma_i)_{i \in \{p, \dots, n\}}$ be the subsequence of $\sigma$ starting at position $p$, so that $\alpha \cdot \beta = \sigma$.
    
    Since $\sigma$ avoids the pattern 010, $\alpha$ and $\beta$ avoid the pattern 010 and they do not share any common values. In particular, $\alpha \in \mathfrak{B}_{p-1,j,i}$ for some $j < m$ and $i < d$. Let us now assume we have already chosen $n,m,d,p,j,i$ and $\alpha$, and count how many sequences $\beta$ fit with those choices.
    
    The subsequence $\beta$ is a 010-avoiding word of length $n-p+1$, which contains exactly $d-i$ distinct values chosen from the remaining $m+1-i$ (that is, all values $\{0, \dots, m\}$ except for the $i$ values in $\alpha$), and such that $\beta_1 = m = \max(\beta)$. Since $m$ is always in $\beta$, it remains to choose the other $d-i-1$ values, there are therefore $\binom{m-i}{d-i-1}$ possible choices for the set of values of $\beta$.
    Once the values of $\beta$ are chosen, there are $\Ustirling{n-p+1}{n-p-d+i+2}$ ways to arrange them into a word avoiding 010 and starting with its maximum, according to Lemma \ref{lem010}. The recursive formula \eqref{eq010} is then obtained by summing over all possible values of $p, i$, and $j$ in this decomposition.
\end{proof}
Using Equation \eqref{eq010}, we generate the terms $\mathfrak{b}_{n,m,d}$ for $n \leqslant 140$ in 1 minute with a C++ program running on a personal computer. Summing over all $m$ and $d$ gives the number of inversion sequences of length $n$ avoiding 010. The first 14 terms of this sequence are
$$1, 2, 5, 15, 53, 215, 979, 4922, 26992, 159958, 1016784, 6890723, 49534501, 376081602.$$

\section{The pair of patterns \{010,000\}} \label{000}
We now enumerate inversion sequences avoiding both patterns 010 and 000. The addition of the pattern 000 does not change our decomposition, so this section is similar to Section \ref{010}. We begin by enumerating a family of words avoiding $\{010,000\}$, like we did for the pattern 010 in Section \ref{010}.

\begin{lem} \label{lem000}
Let $\mathfrak{C}_{n,k}$ be the set of $\{010, 000\}$-avoiding words $\omega = (\omega_1, \dots, \omega_n)$ of length $n$ on the alphabet $\{1, \dots, k\}$ such that $\omega$ contains all letters $\{1, \dots, k\}$ and $\omega_1 = k$. Let $\mathfrak{c}_{n,k} = \#\mathfrak{C}_{n,k}$. The sequence $\mathfrak{c}_{n,k}$ satisfies the following recurrence relation for all $n, k \geqslant 2$:
\begin{equation}
    \mathfrak{c}_{n,k} = (n-1) \cdot \mathfrak{c}_{n-1,k-1} + (n-2) \cdot \mathfrak{c}_{n-2,k-1}.
    \label{eqlem000}
\end{equation}
\end{lem}
\begin{rem}
    Because the words in $\mathfrak{C}_{n,k}$ avoid the pattern 000, we have $\mathfrak{c}_{n,k} = 0$ if $n > 2k$.
\end{rem}
\begin{proof}[Proof of Lemma \ref{lem000}]
    Let $n \geqslant k \geqslant 2$, and $\omega \in \mathfrak{C}_{n,k}$. Since $\omega$ avoids $010$, all letters 1 in $\omega$ are consecutive. Since $\omega$ avoids $000$, $\omega$ has at most two letters 1.
    \begin{itemize}
    \item If $\omega$ contains a single letter 1, then removing it and subtracting $1$ from all other letters yields a word $\omega' \in \mathfrak{C}_{n-1,k-1}$. There are $n-1$ possible positions for the letter 1 in $\omega$, so exactly $n-1$ words $\omega \in \mathfrak{C}_{n,k}$ yield the same $\omega'$.
    \item If $\omega$ contains two letters 1, then removing them and subtracting $1$ from all other letters yields a word $\omega' \in \mathfrak{C}_{n-2,k-1}$. There are $n-2$ possible positions for the two consecutive letters 1 in $\omega$, so exactly $n-2$ words $\omega \in \mathfrak{C}_{n,k}$ yield the same $\omega'$.
    \end{itemize}
    Hence the recurrence relation \eqref{eqlem000} holds for all $n \geqslant k \geqslant 2$.
\end{proof}

\begin{thm} \label{thm000}
Let $\mathfrak{D}_{n,m,d}$ be the set of $\{010, 000\}$-avoiding inversion sequences that have length $n$, maximum $m$, and contain exactly $d$ distinct values. Let $\mathfrak{d}_{n,m,d} = \#\mathfrak{D}_{n,m,d}$. The numbers $\mathfrak{d}_{n,m,d}$ satisfy the following recurrence relation:
\begin{equation}
\mathfrak{d}_{n,m,d} = \sum_{i = 0}^{d-1} \binom{m-i}{d-i-1} \sum_{p = m+1}^n \mathfrak{c}_{n-p+1,d-i} \sum_{j = 0}^{m-1} \mathfrak{d}_{p-1,j,i}.
\label{eq000}
\end{equation}
\end{thm}
\noindent The proof of Theorem \ref{thm000} is identical to that of Theorem \ref{thm010}: the pattern 000 could not spread over $\alpha$ and $\beta$ since they have no value in common.

\section{The pair of patterns \{010,120\}} \label{120}
In this section, we will enumerate inversion sequences avoiding the patterns 010 and 120. We will first solve the enumeration of 'classical' words avoiding this pair of patterns (i.e. words defined by their length and the size of their alphabet). The proof still makes use of words that are required to contain all letters of their alphabet. However, unlike in Sections \ref{010} and \ref{000}, we no longer require the left-most value of our words to be their maximum.
\begin{lem} \label{lem120}
Let $\mathfrak{E}_{n,k}$ be the set of $\{010, 120\}$-avoiding words of length $n$ over the alphabet $\{1, \dots, k\}$. Let $\mathfrak{e}_{n,k} = \#\mathfrak{E}_{n,k}$. Then
\begin{equation}
\mathfrak{e}_{n,k} = \sum_{d=0}^k \binom{k}{d} \frac{\binom{n-1}{d-1} \binom{n+d}{d-1}}{d}.
\label{eqlem120}
\end{equation} 
\end{lem}
\begin{proof}
    Let $\mathfrak{F}_{n,k}$ be the subset of $\mathfrak{E}_{n,k}$ of words that contain all letters $\{1, \dots, k\}$, and $\mathfrak{f}_{n,k} = \#\mathfrak{F}_{n,k}$. By summing over the number of distinct letters $d$ in words in $\mathfrak{E}_{n,k}$, we have the formula
    \begin{equation}
    \mathfrak{e}_{n,k} = \sum_{d=0}^k \binom{k}{d} \mathfrak{f}_{n,d}.
    \label{eqWd120}
    \end{equation}
    We now decompose a word $\omega \in \mathfrak{F}_{n,k}$ around the left-most position of its maximum, similarly to what we did with inversion sequences earlier. Let $p$ be the left-most position of the letter $k$ in $\omega$. Let $\alpha = (\omega_i)_{i \in \{1, \dots, p-1\}}$, and $\gamma = (\omega_i)_{i \in \{p+1, \dots, n\}}$ be two subwords of $\omega$, so that $\alpha \cdot k \cdot \gamma = \omega$.
    
    Since $\omega$ avoids 010 and 120, $\alpha$ has smaller letters than $\gamma$. If $q$ is the maximum of $\alpha$, then the letters of $\alpha$ are exactly $\{1, \dots, q\}$, and those of $\gamma$ are either $\{q+1, \dots, k\}$ or $\{q+1, \dots, k-1\}$. We obtain the recurrence relation
    \begin{equation}
    \mathfrak{f}_{n,k} = \sum_{p=1}^{n} \sum_{q=0}^{k-1} \mathfrak{f}_{p-1,q} \cdot (\mathfrak{f}_{n-p,k-q} + \mathfrak{f}_{n-p,k-q-1}).
    \label{convo120}
    \end{equation}
    Let
    $$F(x,y) = \sum_{n,k \geqslant 0} x^n y^k \mathfrak{f}_{n,k}$$
    be the bivariate ordinary generating function of $\mathfrak{f}_{n,k}$. The recurrence relation \eqref{convo120} translates into the following functional equation for $F$:
    \begin{equation}
    (x + xy) \cdot F(x,y)^2 - (x+1) \cdot F(x,y) + 1 = 0.
    \label{funeq120}
    \end{equation}
    From \eqref{funeq120}, we observe that the function $x^2 \cdot F(x,y)$ satisfies a functional equation for the ordinary generating function of the number of diagonal dissections of a convex $n$-gon into $k$ regions, found in \cite[Section 3.1]{Flajolet_Noy_1999}. Hence $\mathfrak{f}_{n,k}$ is the number of diagonal dissections of a convex $(n+2)$-gon into $k$ regions. A closed formula for the number of diagonal dissections is also obtained in \cite[Section 3.1]{Flajolet_Noy_1999}, which gives
    \begin{equation}
    \mathfrak{f}_{n,k} = \frac{\binom{n-1}{k-1} \binom{n+k}{k-1}}{k}.
    \label{eqbinom120}
    \end{equation}
    Finally, we obtain \eqref{eqlem120} from \eqref{eqWd120} and \eqref{eqbinom120}.
\end{proof}

The enumeration of inversion sequences avoiding 010 and 120 does not require us to refine sequences according to their number of distinct values. In this simpler situation, both our decomposition and the resulting formula are similar to what was done by Mansour and Shattuck in \cite{Mansour_Shattuck_2015} for the single pattern 120.
\begin{thm} \label{thm120}
Let $\mathfrak{G}_{n,m}$ be the set of $\{010, 120\}$-avoiding inversion sequences that have length $n$ and maximum $m$. Let $\mathfrak{g}_{n,m} = \#\mathfrak{G}_{n,m}$. The numbers $\mathfrak{g}_{n,m}$ satisfy the following recurrence relation:
\begin{equation}
\mathfrak{g}_{n,m} = \sum_{p = m+1}^n \sum_{j = 0}^{m-1} \mathfrak{g}_{p-1,j} \cdot \mathfrak{e}_{n-p, m-j}.
\label{eq120}
\end{equation}
\end{thm}
\begin{proof}
    We proceed similarly to the proof of Theorem \ref{thm010}. Let $\sigma \in \mathfrak{G}_{n,m}$. Let $p$ be the left-most position of the value $m$ in $\sigma$. Let $\alpha = (\sigma_i)_{i \in \{1, \dots, p-1\}}$ be the subsequence of $\sigma$ to the left of the left-most $m$, and $\gamma = (\sigma_i)_{i \in \{p+1, \dots, n\}}$ be the subsequence of $\sigma$ to the right of the left-most $m$, so that $\alpha \cdot m \cdot \gamma = \sigma$.

    Since $\sigma$ avoids the patterns 010 and 120, the subsequences $\alpha$ and $\gamma$ avoid these two patterns as well, and all values of $\gamma$ are greater than all values of $\alpha$. We have $\alpha \in \mathfrak{G}_{p-1,j}$ for some $j < m$, and $\gamma$ is a $\{010, 120\}$-avoiding word of length $n-p$ over the alphabet $\{j+1, \dots, m\}$.
    The recursive formula \eqref{eq120} is obtained by summing over all possible $p, j, \alpha$, and $\gamma$.
\end{proof}

\section{The pairs \{010,201\} and \{010,210\}} \label{210}

\indent Wilf-equivalence between the pairs of patterns $\{010,201\}$ and $\{010,210\}$ was proved by Yan and Lin in \cite{Yan_Lin_2020} through the use of a bijection. In this section we will work with the pair of patterns $\{010,210\}$, although our method can also be applied to the pair $\{010,201\}$, resulting in the same formula.

Earlier, we refined the enumeration of pattern-avoiding inversion sequences $\sigma$ according to two parameters (in addition to their size). The first one is the maximum $m$ of $\sigma$, and is needed in order to decompose sequences at the left-most position of their maximum. The second one records the number of values that are forbidden to the right of the left-most maximum, provided $\sigma$ is the left part in our usual decomposition around the left-most occurrence of the maximum (denoted $\alpha$ earlier).

If $\sigma$ is a sequence avoiding a set of patterns $P$, we say that a value $v \in \{0, \dots, \max(\sigma)\}$ is \emph{forbidden} by $\sigma$ and $P$ if $\sigma \cdot m v$ contains a pattern in $P$ when $m > \max(\sigma)$. We denote by $\forb(\sigma, P)$ the number of values forbidden by $\sigma$ and $P$, or simply $\forb(\sigma)$ when there is no ambiguity.

For example, for the patterns 010 and $\{010,000\}$, all values occurring in $\sigma$ were forbidden, so $\forb(\sigma)$ was simply the number of distinct values in $\sigma$, which was recorded by the parameter $d$ in Theorems \ref{thm010} and \ref{thm000}. For the pair $\{010,120\}$, all values from 0 to $\max(\sigma)$ were forbidden, so the number $\max(\sigma)+1$ of forbidden values was redundant with recording the maximum of $\sigma$, and therefore omitted.

\begin{rem} \label{rem210}
    If $\sigma$ is a $\{010, 210\}$-avoiding sequence with values in $\mathbb N$, then $\forb(\sigma) = q+r$, where $q$ is the largest value of $\sigma$ such that a larger value appears to its left, or $q = 0$ if there is no such value (i.e. if $\sigma$ is nondecreasing), and $r$ is the number of values in $\sigma$ that are greater than or equal to $q$ (counted without multiplicity).
\end{rem}
\begin{proof}
All values smaller than $q$ are forbidden by the pattern 210 (regardless of whether or not they occur in $\sigma$). In addition, the value $q$ itself and any values in $\sigma$ greater than $q$ (counted by $r$) are forbidden by the pattern 010.
\end{proof}

\begin{lem} \label{lem210}
Let $\mathfrak{H}_{n,k}$ be the set of $010$-avoiding words $\omega$ of length $n$ over the alphabet $\{1, \dots, k\} \sqcup \{\infty\}$ (where $\infty$ is the largest letter) such that the subword $\omega'$ defined by removing all letters $\infty$ from $\omega$ is nondecreasing, and $k$ is the maximum of $\omega'$ if $\omega'$ is nonempty, or $k = 0$ otherwise.
Consider the bipartition of $\mathfrak{H}_{n,k}$ into $\mathfrak{H}^{(1)}_{n,k} \sqcup \mathfrak{H}^{(2)}_{n,k}$, where $\mathfrak{H}^{(1)}_{n,k} = \{\omega \in \mathfrak{H}_{n,k} \; | \; \omega_n = \infty\}$, and $\mathfrak{H}^{(2)}_{n,k} = \{\omega \in \mathfrak{H}_{n,k} \; | \; \omega_n \neq \infty\}$.
Let $\mathfrak{h}_{n,k} = \#\mathfrak{H}_{n,k}$, $\mathfrak{h}^{(1)}_{n,k} = \#\mathfrak{H}^{(1)}_{n,k}$, and $\mathfrak{h}^{(2)}_{n,k} = \#\mathfrak{H}^{(2)}_{n,k}$. These numbers satisfy the following equations:
\begin{equation}
\begin{cases}
\mathfrak{h}_{n,k} = \mathfrak{h}^{(1)}_{n,k} + \mathfrak{h}^{(2)}_{n,k} \\[5pt]
\mathfrak{h}^{(1)}_{n,k} = \mathfrak{h}_{n-1,k} \\[5pt]
\mathfrak{h}^{(2)}_{n,k} = \mathfrak{h}^{(2)}_{n-1,k} + \sum_{i=0}^{k-1} \mathfrak{h}_{n-1,i}.
\end{cases}
\label{eqlem210}
\end{equation}
\end{lem}
\begin{proof}
    By definition, $\mathfrak{H}_{n,k} = \mathfrak{H}^{(1)}_{n,k} \sqcup \mathfrak{H}^{(2)}_{n,k}$, therefore $\mathfrak{h}_{n,k} = \mathfrak{h}^{(1)}_{n,k} + \mathfrak{h}^{(2)}_{n,k}$.
    Let $\omega \in \mathfrak{H}_{n,k}$.
    \begin{itemize}
        \item If $\omega \in \mathfrak{H}^{(1)}_{n,k}$, then $\omega_n = \infty$. Removing $\omega_n$ yields a word in $\mathfrak{H}_{n-1,k}$.
        \item If $\omega \in \mathfrak{H}^{(2)}_{n,k}$, then by the nondecreasing property, $\omega_n = k$. 
        \begin{itemize}
            \item If $\omega_{n-1} = k$, removing $\omega_n$ yields a word in $\mathfrak{H}^{(2)}_{n-1,k}$.
            \item If $\omega_{n-1} = \infty$, the avoidance of the pattern 010 ensures $\omega$ cannot contain another letter $k$, therefore removing $\omega_n$ yields a word in $\mathfrak{H}^{(1)}_{n-1,i}$ for some $i < k$.
            \item Otherwise, $\omega_{n-1} = i$ for some $i < k$, and removing $\omega_n$ yields a word in $\mathfrak{H}^{(2)}_{n-1,i}$.
        \end{itemize}
    \end{itemize}
    All maps above are bijections, so the last two formulas in \eqref{eqlem210} naturally follow, observing that $\mathfrak{h}^{(1)}_{n-1,i} + \mathfrak{h}^{(2)}_{n-1,i} = \mathfrak{h}_{n-1,i}$.
\end{proof}

\begin{thm} \label{thm210} % corriger ça et changer le programme aussi
Let $\mathfrak{I}_{n,m,f}$ be the set of $\{010, 210\}$-avoiding inversion sequences $\sigma$ that have length $n$, maximum $m$, and such that $\forb(\sigma) = f$. Let $\mathfrak{i}_{n,m,f} = \#\mathfrak{I}_{n,m,f}$.
The numbers $\mathfrak{i}_{n,m,f}$ satisfy the following recurrence relation:
\begin{equation}
\mathfrak{i}_{n,m,f} = \sum_{p = m+1}^n \sum_{i=0}^{f-1} \sum_{j = 0}^{m-1}  \mathfrak{i}_{p-1,j,i} \cdot \mathfrak{h}_{n-p,f-i-1}.
\label{eq210}
\end{equation}
\end{thm}
\begin{proof}
    We proceed the same way as before. Let $\sigma \in \mathfrak{I}_{n,m,f}$. Let $p$ be the left-most position of the value $m$ in $\sigma$. Let $\alpha = (\sigma_i)_{i \in \{1, \dots, p-1\}}$ be the subsequence of $\sigma$ to the left of the left-most $m$, and $\gamma = (\sigma_i)_{i \in \{p+1, \dots, n\}}$ be the subsequence of $\sigma$ to the right of the left-most $m$, so that $\alpha \cdot m \cdot \gamma = \sigma$. Let $\gamma'$ be the subsequence of $\gamma$ obtained by removing all values $m$ from $\gamma$ (if $\gamma$ does not contain $m$, then $\gamma' = \gamma$).
    
    Since $\sigma$ avoids the pattern 010, the subsequences $\alpha$ and $\gamma$ avoid 010, and they do not share any common values. Since $\sigma$ avoids the pattern 210, $\alpha$ avoids 210, and $\gamma'$ is nondecreasing.

    The subsequence $\alpha$ is in $\mathfrak{I}_{p-1,j,i}$ for some $j < m$ and $i \leqslant f$. Notice that we actually have $i < f$ since $m$ is a forbidden value for $\sigma$ (because of the avoidance of 010), but not for $\alpha$ (since $m > \max(\alpha))$. Let $s$ be the largest value of $\alpha$ such that a larger value appears to its left, or $s = 0$ if there is no such value. Let $t$ be the number of values in $\alpha$ that are greater than or equal to $s$, counted without multiplicity.  By Remark \ref{rem210}, we have $s+t = \forb(\alpha) = i$.
    
    The subsequence $\gamma$ is a 010-avoiding word of length $n-p$ over an alphabet $\Sigma$ of size $m-i+1$ (that is, all values in $\{s, \dots, m\}$ except the $t$ values in $\{s, \dots, m\}$ forbidden by $\alpha$), and such that $\gamma'$ is nondecreasing. Since $f$ values are forbidden by $\sigma$, including $i$ values previously forbidden by $\alpha$, there are $f-i$ new values forbidden by $m \cdot \gamma$.
    The new values forbidden by $m \cdot \gamma$ are precisely:
    \begin{itemize}
        \item Every value in $\Sigma$ lower than the maximum of $\gamma'$, because of the avoidance of 210 (with $m$ playing the role of the 2, and the maximum of $\gamma'$ playing the role of the 1).
        \item The value $m$, and the maximum of $\gamma'$, because of the avoidance of 010.
    \end{itemize}
    It follows that the maximum of $\gamma'$ is the $(f-i-1)$-th smallest letter of $\Sigma$, therefore the number of possible choices for $\gamma$ is counted by $\mathfrak{h}_{n-p,f-i-1}$ (with $m$ taking the role of the letter denoted $\infty$ in Lemma \ref{lem210}).
    As before, the recursive formula \eqref{eq210} is obtained by summing over all possible $p, i, j, \alpha$, and $\gamma$.
\end{proof}

\section{The pair of patterns \{010,110\}} \label{110}
In the case of the pair $\{010,110\}$, we will apply our decomposition of 010-avoiding inversion sequences to words as well. This can also be done with the words studied in Sections \ref{010}, \ref{000}, and \ref{120}, however the methods used in those sections resulted in simpler, more efficient formulas. We begin with a characterization of forbidden values similar to Remark \ref{rem210}.

\begin{rem} \label{rem110}
    If $\sigma$ is a $\{010, 110\}$-avoiding sequence with values in $\mathbb N$, then $\forb(\sigma) = q+r$, where $q$ is the largest value which appears twice in $\sigma$, or $q = 0$ if there is no such value, and $r$ is the number of values in $\sigma$ that are greater than or equal to $q$ (counted without multiplicity).
\end{rem}
\begin{proof}
All values smaller than $q$ are forbidden by the pattern 110 (regardless of whether or not they occur in $\sigma$). In addition, the value $q$ itself and any values in $\sigma$ greater than $q$ (counted by $r$) are forbidden by the pattern 010.
\end{proof}

\noindent Let $\delta_{a,b} = \left \{ \begin{array}{lcl}
1 & \text{if} & a = b \\
0 & \text{if} & a \neq b
\end{array} \right.$ be the Kronecker delta function.
\begin{lem} \label{lem110}
Let $\mathfrak{J}_{n,k,f}$ be the set of $\{010,110\}$-avoiding words $\omega$ of length $n$ over the alphabet\footnote{Our alphabet now starts at 0 so that our characterization of forbidden values in Remark \ref{rem110} remains the same for both words and inversion sequences.} $\{0, \dots, k-1\}$ such that $\forb(\omega) = f$. Let $\mathfrak{K}_{n,k}$ be the set of $\{010,110\}$-avoiding words $\omega$ of length $n$ over the alphabet $\{0, \dots, k-1\}$.
Let $\mathfrak{j}_{n,k,f} = \#\mathfrak{J}_{n,k,f}$, and $\mathfrak{k}_{n,k} = \#\mathfrak{K}_{n,k}$, so that we have
$$\mathfrak{k}_{n,k} = \sum_{f=0}^k \mathfrak{j}_{n,k,f}.$$
The numbers $\mathfrak{j}_{n,k,f}$ satisfy the following recurrence relation:
\begin{equation}
\mathfrak{j}_{n,k,f} = \sum_{p=1}^{n} \sum_{i=0}^{f-1} \sum_{m=0}^{k-1} \mathfrak{j}_{p-1,m,i} \cdot (\mathfrak{j}_{n-p, m-i, f-i-1} + \delta_{f, m+1} \cdot \sum_{\ell=0}^{n-p-1} \mathfrak{k}_{\ell, m-i}).
\label{eqlem110}
\end{equation}
\end{lem}
\begin{proof}
    Let $\omega \in \mathfrak{J}_{n,k,f}$. Let $m$ be the maximum of $\omega$ and $p$ the position of the left-most $m$ in $\omega$. Let $\alpha = (\omega_i)_{i \in \{1, \dots, p-1\}}$, and $\gamma = (\omega_i)_{i \in \{p+1, \dots, n\}}$, so that $\alpha \cdot m \cdot \gamma = \omega$.

    Since $\omega$ avoids the pattern 010, the subwords $\alpha$ and $\gamma$ avoid 010, and they do not share any common letters. Since $\sigma$ avoids the pattern 110, $\alpha$ and $\gamma$ avoid 110, and the letters of $\gamma$ are greater than any letters which appear twice in $\alpha$. Additionally, if a letter $m$ appears in $\gamma$, then all letters to its right must be $m$.

    If $i$ letters are forbidden by $\alpha$, i.e. $\alpha \in \mathfrak{J}_{p-1,m,i}$ then $\gamma$ is a $\{010,110\}$-avoiding word of length $n-p$ on an alphabet of size $m+1-i$ (that is, all letters $\{0, \dots, m\}$ except the $i$ letters forbidden by $\alpha$), and such that all letters to the right of a $m$ are also $m$. We distinguish two cases for $\gamma$:
    \begin{itemize}
        \item If $\gamma$ does not contain the letter $m$, then $\gamma$ is a $\{010,110\}$-avoiding word of length $n-p$ on an alphabet of size $m-i$. Since $f$ letters are forbidden by $\omega$, including $i$ letters forbidden by $\alpha$ and 
        the letter $m$, the remaining $f-i-1$ must be forbidden by $\gamma$, therefore the number of possible choices for $\gamma$ is $\mathfrak{j}_{n-p, m-i, f-i-1}$.
        \item If $\gamma$ contains the letter $m$, let $\ell$ be the number of letters different from $m$ in $\gamma$ (counted with multiplicity). No letters other than $m$ may appear to the right of a $m$ in $\gamma$, therefore $\gamma$ can be seen as a $\{010,110\}$-avoiding word of length $\ell$ on an alphabet of size $m-i$ to which we append $n-p-\ell$ letters $m$. In this case, $m$ is the largest letter which appears twice in $\omega$, and $\omega$ does not contain any letters greater than $m$, therefore $f = m+1$.
    \end{itemize}
    The recursive formula \eqref{eqlem110} is obtained by summing over all $p, i, m, \ell, \alpha$, and $\gamma$.
\end{proof}

\begin{thm} \label{thm110}
Let $\mathfrak{L}_{n,m,f}$ be the set of $\{010,110\}$-avoiding inversion sequences $\sigma$ of length $n$, maximum $m$, and such that $\forb(\sigma) = f$. Let $\mathfrak{l}_{n,m,f} = \#\mathfrak{L}_{n,m,f}$.
The numbers $\mathfrak{l}_{n,k,f}$ satisfy the following recurrence relation:
\begin{equation}
\mathfrak{l}_{n,m,f} = \sum_{p=m+1}^{n} \sum_{i=0}^{f-1} \sum_{j=0}^{m-1} \mathfrak{l}_{p-1,j,i} \cdot (\mathfrak{j}_{n-p, m-i, f-i-1} + \delta_{f,m+1} \cdot \sum_{\ell=0}^{n-p-1} \mathfrak{k}_{\ell, m-i}).
\label{eq110}
\end{equation}
\end{thm}
\noindent The proof of Theorem \ref{thm110} is identical to that of Lemma \ref{lem110}, except the left part of our decomposition is now an inversion sequence instead of a word.

\section{The pair of patterns \{010,102\}} \label{102}
The enumeration of inversion sequences avoiding the pair of patterns $\{010, 102\}$ remains unsolved, as our usual method does not apply to it. Indeed, the value 1 in the pattern 102 may appear in the left part of our decomposition, and the values 02 in the right part, hence recording the number of forbidden values in not sufficient in this case.

\section*{Acknowledgements}
The author would like to thank Mathilde Bouvel and Emmanuel Jeandel for their helpful suggestions and advice.

\emergencystretch=1em
\printbibliography

\end{document}